\ifx\documentclass\undefined
\documentstyle[12pt]{article}
\else
\documentclass[12pt]{article}
\usepackage{latexsym}
\usepackage{amsfonts}
\usepackage{amsmath}
\usepackage{amssymb}
\fi

\author{K. Bezdek \thanks{Partially supported by the Hung.
Nat. Sci. Found (OTKA), grant no. T043556} 
\thanks{Partially supported by a Natural Sciences and 
Engineering Research Council of Canada Discovery Grant}
\and A. E. Litvak$\mbox{}^{\ddagger}$ }

\date{}

\font\tenBbb=msbm10 at 12pt         \font\sevenBbb=msbm9    \font\fiveBbb=msbm7
\newfam\Bbbfam
\textfont\Bbbfam=\tenBbb \scriptfont\Bbbfam=\sevenBbb
\scriptscriptfont\Bbbfam=\fiveBbb

\def\R{{\mathbb R}}
\def\S{{\mathbb S}}

\def\no{ \nor \cdot  \nor }

\def\kkk{\null\hfill $\Box$\smallskip}
\def\r{ \right}

\def\eps{\varepsilon}

\def\lam{\lambda}

\def\nor{ \| }

\def\la{\left\langle}
\def\ra{\r\rangle}

\def\KK{{\bf K}}
\def\LL{{\bf L}}

\def\PP{{\bf P}}
\def\BB{{\bf B}}
\def\CC{{\bf C}}

\def\TT{{\bf T}}

\newcommand{\proof}{{\noindent\bf Proof:{\ \ }}}

\newtheorem{theorem}{Theorem}[section]
\newtheorem{lemma}[theorem]{Lemma}

\newtheorem{prop}[theorem]{Proposition}
\newtheorem{sled}[theorem]{Corollary}
\newtheorem{con}[theorem]{Conjecture}

\newcommand{\conv}{\mathop{\rm conv\,}}
\newcommand{\absconv}{\mathop{\rm abs \,conv\,}}
\newcommand{\vol}{\mathop{\rm vol\,}}

\title{On the vertex index of convex bodies 
\footnote{Keywords: convex body, illumination parameter, vertex index, 
Boltyanski-Hadwiger conjecture, volume ratio.  
2000 Mathematical Subject Classification. Primary: 46B07, 46B09, 52A.
Secondary: 51M16, 53A55}}

\begin{document}

\maketitle

\begin{abstract}
We introduce the vertex index, ${\rm vein} (\KK)$, of a given centrally symmetric 
convex body $\KK \subset \R^d$, which, in a sense, measures how well $\KK$ can be 
inscribed into a convex polytope with small number of vertices. This index is closely 
connected to the illumination parameter of a body, introduced earlier by the first named 
author, and, thus, related to the famous conjecture in Convex Geometry about covering 
of a $d$-dimensional body by $2^d$ smaller positively homothetic copies. We 
provide asymptotically sharp estimates (up to a logarithmic term) of this index 
in the general case. More precisely, we show that for every centrally symmetric 
convex body $\KK \subset \R^d$ one has 
$$
  \frac{  d^{3/2}}{\sqrt{2 \pi e}\ 
  \mbox{\rm ovr} (\KK) } \leq {\rm vein}(\KK)\leq C\ d^{3/2}  \ \ln(2d) ,
$$
where $\mbox{ovr} (\KK) = \inf \left( \vol ({\cal E}) / \vol (\KK) \r)^{1/d}$
is the outer volume ratio of $\KK$ with the infimum taken over all ellipsoids 
${\cal E} \supset \KK$ and with $\vol (\cdot)$ denoting the volume. 
Also, we provide sharp estimates in dimensions 2 and 3. Namely, in 
the planar case we prove that $4  \leq {\rm vein}(\KK)\leq 6$ 
with  equalities for parallelograms and 
affine regular convex hexagons, 
and in the $3$-dimensional case we show that 
$6  \leq {\rm vein}(\KK)$  with equality for octahedra. 
We conjecture that the vertex index of a $d$-dimensional Euclidean ball 
(resp., ellipsoid) is $2 d\sqrt{d}$. We prove this  conjecture in dimensions 
two and three.
\end{abstract}

\section{Introduction}
\label{zero}
Let $\KK$ be a convex body symmetric about the origin $0$ in $\R^d, d\ge 2$
(such bodies below we call $0$-symmetric convex bodies). Now, we place $\KK$ 
in a convex polytope, say $\PP$, with vertices $p_1, p_2, \dots , p_n$, where 
$n\ge d+1$. Then it is natural to measure the closeness of the vertex set of 
$\PP$ to the origin $0$ by computing $\sum_{1\le i\le n}{\| p_i\|}_{\KK }$, 
where ${\| x\|}_{\KK}=\inf \{\lambda >0 \ |\ x\in\lambda\KK\}$ denotes the norm 
of $x\in \R^d$ generated by $\KK$. Finally, we look for the convex polytope that 
contains $\KK$ and whose vertex set has the smallest possible closeness to $0$ 
and introduce the {\it vertex index}, ${\rm vein}(\KK)$, of $\KK$ as follows:
$$
   {\rm vein}(\KK )=\inf \left\{ \sum_{i}{\| p_i\|}_{\KK } \ |\ 
   \KK\subset \conv\{p_i\} \right\}.
$$

We note that
${\rm vein}(\KK )$ is an affine invariant quantity assigned to $\KK$, 
i.e. if $A : \R^d \to \R^d$ is an (invertible) linear map, then 
${\rm vein}(\KK )={\rm vein}(A(\KK) )$.
The main goal of this paper is to give lower and upper estimates on ${\rm vein}(\KK )$. 
This question seems to raise a fundamental problem that is connected to some important 
problems of analysis and geometry including the problem of estimating the illumination 
parameters of convex bodies, the Boltyanski-Hadwiger illumination conjecture, 
some of the problems on covering a convex body by another one, and the problem 
of estimating the Banach-Mazur distances between convex bodies. 
Section 3 of this paper provides more details on these connections. 
Next we summarize the major results of our paper.

\medskip

\noindent
{\bf Theorem A \ } {\it 
For every $d\geq 2$ one has 
$$
  \frac{  d^{3/2}}{\sqrt{2 \pi e}} 
  \leq {\rm vein}(\BB_2^d)  \leq 2 d^{3/2} ,
$$
where $\BB_2^d$ denotes the Euclidean unit ball in $\R^d$. 
Moreover, if $d=2$, $3$ then ${\rm vein}(\BB_2^d)=2 d^{3/2}$.
}

\medskip

In fact, the above theorem is a combination of Theorem~\ref{Eball} and of 
Corollary~\ref{ballest} in Sections 4 and 5. In connection with that it seems 
natural to conjecture the following.

\medskip

\noindent
{\bf Conjecture B\ } {\it 
For every $d\geq 2$ one has 
$${\rm vein}(\BB_2^d)=2 d^{3/2}.$$
}

If Conjecture~B holds, then it is easy to see that it implies via 
Lemma~\ref{distest} the inequality ${\rm vein}(\KK ) \ge 2d$ for 
any $0$-symmetric convex body $\KK$ in $\R^d$. This estimate was recently 
obtained in \cite{GL}. Note that by Proposition~\ref{octah} below, 
${\rm vein}(\CC)=2d$, where $\CC$ denotes any $d$-dimensional 
crosspolytope of $\R^d$.

The following is the major result of Section 5, which is, in fact, a combination of 
Theorems~\ref{ball} and~\ref{ru}. 

\medskip 

\noindent
{\bf Theorem C \ }{\it 
There are absolute constants $c>0, C>0$ such that for every 
$d\ge 2$ and every $0$-symmetric convex body $\KK$  in $\R^d$ one has
$$
  \frac{  d^{3/2}}{\sqrt{2 \pi e}\ \mbox{\rm ovr} (\KK) }   \leq 
   {\rm vein}(\KK)\leq C\ d^{3/2}  \ \ln(2d),
$$
where $\mbox{\rm ovr} (\KK) = \inf \left( \vol ({\cal E}) / \vol (\KK) \r)^{1/d}$ 
is the outer volume ratio of $\KK$ with the infimum taken over all ellipsoids  
${\cal E} \supset \KK$ and with $\vol (\cdot)$ denoting the volume. 
}

\medskip

Examples of a cross-polytope $\CC$ (see Proposition~\ref{octah}) and of $\BB _2^d$ 
(see Theorem~A) show that both estimates in Theorem~C can be asymptotically sharp, 
up to a logarithmic term. One may wonder about the {\it precise} bounds. 
Section 4 investigates this question in dimensions $2$ and $3$. However, in high 
dimensions the answer to this question might be different. As we mentioned above, 
the function ${\rm vein} (\cdot )$ 
attains its minimum at crosspolytopes. It is not clear to us for what convex bodies 
should the function ${\rm vein} (\cdot )$ attain its maximum. In particular, as 
Corollary~\ref{ballest} gives an upper estimate on the vertex index of $d$-cubes 
which is somewhat weaker than the similar estimate for Euclidean $d$-balls, it is 
natural to ask, whether the function ${\rm vein} (\cdot )$ attains its maximum at 
(affine) cubes (at least in some dimensions). On the other hand, it would not come 
as a surprise to us if the answer to this question were negative, in which case it 
seems reasonable to suggest the ellipsoids (in particular, in dimensions of the form 
$d=2^m$) or perhaps, the dual of $\S - \S$, where $\S$ denotes any simplex, as convex 
bodies for which the function ${\rm vein} (\cdot)$ attains its maximum. 

\medskip

\noindent
{\bf Acknowledgment. } In the first version of this paper the estimate 
$$
  \frac{ c \ d^{3/2}}{\mbox{\rm ovr} (\KK) \ \sqrt{\ln (2 d)}} 
  \leq {\rm vein}(\KK)\leq C\ d^{3/2}  \ \ln(2d)
$$
in Theorem~C was proved. The proof used a volumetric result from \cite{GMP}, 
the Hadamard inequality, and an averaging argument. E.~D.~Gluskin noticed that 
the use of a result from \cite{BP} (together with Santal{\'o} inequality) instead 
allows to remove the logarithmic term in the lower bound. We are grateful to 
Gluskin for this remark. We are also grateful to the anonymous referee, who 
later noticed the same improvement.

\section{Notations}
\label{one}

In this paper we identify the a $d$-dimensional affine space with $\R ^d$. 
By $|\cdot|$ and $\la \cdot , \cdot \ra$ we denote the canonical 
Euclidean norm and the canonical inner product on $\R ^d$.  
The canonical basis of $\R ^d$ we denote by $e_1, \ldots, e_d$. By 
$\| \cdot \| _p$, $1\leq p \leq \infty$, we denote the $\ell _p$-norm, 
i.e. 
$$\|x\| _p = \left( \sum _{i\geq 1} |x_i|^p \r) ^{1/p} \, 
\mbox{ for } \ p < \infty \, \, \, \, \mbox{ and } \, \, \, \, 
\|x\| _{\infty } = \sup _{i\geq 1} |x_i|.
$$  
In particular, $\no _2 = |\cdot |$. 
As usual, $\ell _p^d = (\R ^d, \|\cdot \|_p)$,  and  
the unit ball of $\ell _p^d$ is denoted  by $\BB_p^d$. 

Given points $x_1, \ldots, x_k$ in $\R ^d$ we denote their 
convex hull by $\conv \{x_i\}_{i\leq k}$ and their absolute 
convex hull by $\absconv \{x_i\}_{i\leq k} = 
\conv \{\pm x_i\}_{i\leq k}$.  
Similarly, the convex hull of a set $A\subset\R^d$ is denoted 
by $\conv A$ and absolute  convex hull of $A$ is denoted 
by $\absconv A$ ($=\conv A\cup -A$).  

Let $\KK\subset \R^d$ be a convex body, i.e. a compact convex set 
with non-empty interior such that the origin $0$ of $\R^d$ belongs 
to $\KK$. We denote by $\KK^{\circ}$ the polar of $\KK$, i.e.
$$
   \KK^{\circ} = \left\{ x \, \, | \, \, \la x, y \ra \leq 1 \, \,
   \mbox{ for every } \, \, y \in \KK \r\}.
$$
As is well-known, if $E$ is a linear subspace of $\R^d$, then
the polar of $\KK \cap E$ (within $E$) is 
$$
   \left(\KK \cap E\r) ^{\circ} = P_E \KK ^{\circ},   
$$
where $P_E$ is the orthogonal projection onto $E$. Note also 
that $\KK ^{\circ \circ} = \KK$. 

If $\KK$ is an 0-symmetric convex body, then the Minkowski functional of $K$,  
$$
    {\| x\|}_{\KK}=\inf \{\lambda >0 \ |\ x\in\lambda\KK\}, 
$$ 
defines a norm on $\R^d$ with the unit ball $\KK$. 

The {\it Banach-Mazur distance} between two $0$-symmetric convex  bodies 
$\KK$ and $\LL$ in $\R ^d$ is defined by
$$
     d(\KK, \LL) = \inf{\left\{ \lam > 0       \ \mid \ 
      \LL \subset T \KK \subset        \lam \LL \r\}},
$$
where the infimum is taken over all linear operators 
$T : \R^d \to \R^d$. It is easy to see that 
$$
   d(\KK, \LL) = d(\KK^{\circ}, \LL^{\circ}) . 
$$
The Banach-Mazur distance between $\KK$ and the closed Euclidean 
ball $\BB_2^d$ we denote by $d_{\KK}$. As it is well-known, 
John's Theorem (\cite{J}) implies that for every $0$-symmetric convex body 
$\KK$, $d_{\KK}$ is bounded by $\sqrt{d}$. 
Moreover, $d_{\BB_1^d} = d_{\BB_{\infty}^d} = \sqrt{d}$ 
(see e.g. \cite{T}).

Given a (convex) body $\KK$ in $\R^d$ we denote its volume by $\vol (\KK)$. 
Let $\KK$ be a $0$-symmetric convex body in $\R^d$. 
The {\it outer volume ratio} of $\KK$ is 
$$
 \mbox{ovr} (\KK) = \inf \left( \frac{\vol ({\cal E})}{\vol (\KK) }\r)^{1/d}, 
$$
where the infimum is taken over all  $0$-symmetric ellipsoids in $\R^d$ 
containing $\KK$. By John's theorem we have
$$
   \mbox{ovr} (\KK) \leq \sqrt{d}.
$$
Note also that 
$$
   \vol (\BB_2^d) = \frac{\pi ^{d/2} }{ \Gamma (1+d/2)} \leq 
   \left(\frac{2 \pi e}{d} \r)^{d/2}, 
$$
where $\Gamma (\cdot)$ denotes the Gamma-function. 
 
Given a finite set $A$ we denote its cardinality by $|A|$.

\section{Preliminary results and relations to other problems}
\label{two}

Let $\KK$ be a  $0$-symmetric convex body in $\R^d, d\ge 2$. 
An exterior point $p\in \R^d\setminus \KK$ of $\KK$ illuminates a boundary point 
$q$ of $\KK$ if the half line emanating from $p$ passing through $q$ intersects the 
interior of $\KK$ (after the point $q$). 
Furthermore, a family of exterior 
points of $\KK$, say $\{p_1, p_2, \dots , p_n\} \subset \R^d\setminus \KK$, illuminates 
$\KK$ if each boundary point of $\KK$ is illuminated by at least one of the points 
$p_1, p_2, \dots , p_n$. The points $p_1, p_2, \dots , p_n$ here are called light sources. 
The well-known Boltyanski-Hadwiger conjecture says that every $d$-dimensional 
convex body $\KK$ can be illuminated by $2^d$ points. Clearly,  we need $2^d$ points to 
illuminate any $d$-dimensional affine cube. The Boltyanski-Hadwiger conjecture is equivalent 
to another famous long-standing conjecture in Convex Geometry, which says that every 
$d$-dimensional convex body $\KK$ can be covered by $2^d$ smaller positively homothetic copies of 
$\KK$. Again, the example of a $d$-dimensional affine cube shows that  $2^d$ cannot be improved 
in general. We refer the interested reader to \cite{Be93}, \cite{Be05}, \cite{MS} 
for further information and partial results on these conjectures.

Although computing the smallest number of points illuminating a given body is very important,  
it does not provide any quantitative information on points of illumination. In particular, one 
can take light sources to be very far from the body. To control that, 
the first named author introduced (\cite{Be92}) the {\it illumination parameter},  
${\rm ill}(\KK )$, of $\KK$  as follows:
$$
   {\rm ill}(\KK )=\inf\left\{ \sum_{i}{\| p_i\|}_{\KK } \ \ \ |\ \ \ \{p_i\}_i
   \  \text{illuminates} \  \KK\right\}.
$$
Clearly this insures that far-away light sources are penalized. In \cite{Be92} the following 
theorem was stated with an outline of its proof. (The detailed proof can be found in \cite{Be06}). 

\begin{theorem}\label{hexa}
If $\KK$ is a $0$-symmetric convex domain of $\R^2$, then ${\rm ill}(\KK )\le 6$ with equality 
for any affine regular convex hexagon.
\end{theorem}

In the same paper the problem of finding the higher dimensional analogue of that 
claim was raised as well.

\smallskip

Motivated by the notion of the illumination parameter Swanepoel \cite{S} introduced the 
{\it covering parameter}, ${\rm cov}( \KK)$, of $\KK$ in the following way.
$$
  {\rm cov}( \KK)=\inf \left\{ \sum_{i} (1-\lambda_i)^{-1} \ |\  \KK \subset 
  \bigcup _{i} (\lambda_i\KK + t_i), 0<\lambda_i<1, t_i\in \R^d\right\}.
$$
In this way homothets almost as large as $\KK$ are penalized. Swanepoel \cite{S} proved the 
following inequality.

\begin{theorem}
There exists an absolute constant $C$ such that for every $0$-symmetric convex body 
$\KK$ in $\R^d$, $d\geq 2$ one has  
$$
  {\rm ill}(\KK )\le 2 \cdot{\rm cov}( \KK)\le C 2^d d^2 \ln d .
$$
\end{theorem}

It is not difficult to see that for any convex body $\KK$ in $\R^d, d\ge 2$ one has 
${\rm vein}(\KK )\le {\rm ill}(\KK )$ with equality for all smooth $\KK$. 
Thus, the above two theorems yield the following immediate result.

\begin{sled}\label{ill}
Let $\KK$ be an $0$-symmetric convex body in $\R^d, d\geq 2$. Then 
\item (i) in case of $d=2$ the inequality ${\rm vein}(\KK )\le 6$ holds;
\item (ii) in case of $d\ge 3$ the inequality 
${\rm vein}(\KK )\le C 2^d d^2 \ln d$ stands.
\end{sled}

As we mentioned in the Introduction, the main goal of this paper is to improve 
the above estimates and also to give lower bounds. We note that 
Theorem~A and C essentially improve the previously known estimates on the illumination 
parameter (of smooth convex bodies). Indeed, they immediately imply the following corollary.

\begin{sled} 
For every $d\geq 2$ and every 
$0$-symmetric convex body $\KK\subset \R^d$ one has 
$$
  \frac{  d^{3/2}}{\sqrt{2 \pi e}\ 
  \mbox{\rm ovr} (\KK) }  \leq {\rm ill}(\KK) .
$$
Moreover, if  $\KK$ is smooth, then
$$
   {\rm ill}(\KK)\leq C\ d^{3/2}  \ \ln(2d),  
$$
where $C>0$ is an absolute constant.
\end{sled}

Finally, we mention two results on Banach-Mazur distances, that will be used below.

\begin{lemma} \label{distest} 
Let $\KK$ and $\LL$ be  $0$-symmetric convex bodies in $\R ^d$.
Then
$$
   {\rm vein} (\KK) \leq d\left(\KK, \LL\r)\cdot {\rm vein} (\LL).
$$
\end{lemma}

\proof Let $T$ be a linear operator such that 
$\KK \subset  T \LL \subset \lam \KK$. Let $p_1, p_2, ..., p_n \in \R^d$ 
be such that $\conv \{p_i\}_{1\le i\le n} \supset \LL$.
Then $\conv \{T p_i\}_{1\le i\le n} \supset T\LL \supset \KK$. Since 
$T\LL \subset \lam \KK$,  we also have $\no _{\KK}\leq \lam \no_{T\LL}$.
Therefore, 
$$
 \sum _{1\le i\le n} \|Tp_i \| _{\KK} \leq \lam \sum _{1\le i\le n} \|Tp_i \| _{T\LL} 
 =  \lam \sum _{1\le i\le n} \|p_i \| _{\LL},   
$$
which implies the desired result.
\kkk

\medskip

\noindent
{\bf Remark. } It is known (\cite{A}, see also \cite{T}) that for every 
$2$-dimensional $0$-symmetric convex body $\KK$ one has 
$d(\KK, \BB_{\infty}^2) \leq 3/2$. Since, clearly, 
${\rm vein} (\BB_{\infty}^2) \leq 4$, we immediately obtain 
$$
  {\rm vein} (\KK) \leq d(\KK, \BB_{\infty}^2) \ {\rm vein} (\BB_{\infty}^2) \leq 6, 
$$
 reproving (i) of Corollary~\ref{ill}.

We will also use the following result (Theorem~2 in \cite{GKM}, see also 
Proposition~37.6 in \cite{T}).

\begin{theorem} \label{cubeoc} For every $d\geq 1$ we have 
$$
   d \left( \BB_1^d, \BB_{\infty}^d \r) \leq C \sqrt{d}, 
$$
with $C=1$ if $d=2^m$ for some integer $m$ and $C = \sqrt{2}+1$ in 
the general case.
\end{theorem}

\section{The vertex index in dimensions $2$ and $3$}

In this section we prove the following theorem.

\begin{theorem}\label{Eball}
\item[(i)] 
For the Euclidean balls in $\R^2$ and $\R^3$ we have 
$${\rm vein}(\BB_2^2)=4\sqrt{2}, \quad \quad \quad \quad {\rm vein}(\BB_2^3)=6\sqrt{3}.$$

\item[(ii)] In general, if $\KK \subset \R^2$, $\LL \subset \R^3$  are arbitrary $0$-symmetric 
convex bodies, then 
$$ 4\le {\rm vein}(\KK)\le 6\le{\rm vein}(\LL)\le 18.$$
\end{theorem}

\medskip

\noindent 
{\bf Remarks.}
\newline
{\bf 1.} Clearly, ${\rm vein} (\BB_1^d) \leq 2d$. Thus, Theorem~\ref{Eball} 
implies ${\rm vein} (\BB_1^d) = 2d$ for $d=2$ and $d=3$. Below (Proposition~\ref{octah}) 
we extend this equality to the general case.
\newline
{\bf 2.} By Remark 1, the  lower estimates in $(ii)$ are sharp. Moreover, it is not hard to 
see that the upper estimate $6$ in the planar case is also sharp by taking any affine regular 
convex hexagon (cf. Theorem~\ref{hexa}).  
\newline
{\bf 3.} We do not know the best possible upper estimate in the $3$-dimensional 
case. It seems reasonable to conjecture the following.

\begin{con}
If $\KK$ is an arbitrary $0$-symmetric convex body in $\R^3$, then 
$${\rm vein}(\KK)\le 12$$
with equality for truncated octahedra of the form $\TT-\TT$, where $\TT$ denotes an arbitrary tetrahedron of $\R^3$.
\end{con}

Note that by Lemma~\ref{distest} this conjecture would be true if, for example, one could 
prove that $d(\BB_1^3, \KK)\leq 2$ for every $0$-symmetric $3$-dimensional convex body. 
To the best of our knowledge no estimates are known for $\max _{\KK}  d(\BB_1^3, \KK)$, 
except the trivial bound $3$. Note also that any bound better than 3 will improve the 
estimate 18 in part (ii) of Theorem~\ref{Eball}.

To prove Theorem~\ref{Eball} we need the following Lemma. The lemma can be proved 
using standard analytic approach or tools like MAPLE. We omit the details.

\begin{lemma}\label{func}
Let $f$ be a function of two variables defined by 
$$
  f(x,y)=\tan{\frac{\pi}{y}} \tan{\left(\frac{x+(y-2)\pi}{2y}\r)}.
$$  
Then 
\item[(i)] for every  fixed $0<x_0< 2\pi$ the function $f(x_0, y)$ 
is decreasing in $y$ over the interval $[3, \infty)$;

\item[(ii)] for every fixed $y_0\ge 3$ the function  
$f(x, y_0)$ is increasing in $x$ over the interval $(0, 2\pi )$;

\item[(iii)] for every fixed $y_0\ge 3$ the function  
$f(x, y_0)$ is convex on the interval $(0, 2\pi )$;

\item[(iv)] $f$ is convex on the  closed rectangle 
$\{(x,y)\ |\ 0.4\le x\le 5.5, 3\le y\le 9\}$.

\end{lemma}

\bigskip

\noindent
{\bf Proof of Theorem~\ref{Eball}:}

$(i)$ The upper estimate ${\rm vein}(\BB_2^d) \leq 2 d \sqrt{d}$ is trivial, since 
$\BB_2^d \subset \sqrt{d} \ \BB_1^d$ for every $d$ (cf. Corollary~\ref{ballest} below). 
We show the lower estimates.

\smallskip 

Let $\PP\subset\R^2$ be a convex polygon with vertices $p_1, p_2, \dots , p_n, n\ge 3$ 
containing $\BB_2^2$. Let $\PP^{\circ}$ denote the polar of $\PP$. Assume that the side 
of $\PP^{\circ}$ corresponding to the vertex $p_i$ of $\PP$ generates the central angle 
$2\alpha_i$ with vertex $0$. Clearly, $0<\alpha_i<{\pi/2}$ and $|p_i| \geq 
\frac{1}{\cos\alpha_i}$ for all $i\le n$. As 
$\frac{1}{\cos x}$ is a convex function over the open interval $(-{\pi/2},{\pi/2})$ 
therefore the Jensen inequality implies that
$$ 
  \sum_{i=1}^{n}|p_i| \geq \sum_{i=1}^{n}\frac{1}{\cos\alpha_i}\ge \frac{n}{\cos \left( 
  \frac{\sum_{i=1}^{n}\alpha_i}{n}\right)}= \frac{n}{\cos{\frac{\pi}{n}}}.
$$
It is easy to see that 
$\frac{n}{\cos{(\pi/n)}}\ge \frac{4}{\cos{(\pi/4)}}=4\sqrt{2}$ 
holds for all $n\ge 3$. 
Thus, ${\rm vein}(\BB_2^2)\ge 4\sqrt{2}$. 
This completes the proof in the planar case.

\medskip

Now, we handle the $3$-dimensional case. Let $\PP\subset\R^3$ be a convex 
polyhedron with vertices $p_1, p_2, \dots , p_n$, $n\ge 4,$ containing $\BB_2^3$. 
Of course, we assume that $|p_i| > 1$. We distinguish the following three cases: 
(a) $n=4$ , (b) $n\ge 8$ and (c) $5\le n\le 7$.
In fact, the proof given for Case (c) works also for Case (b), 
however the Case (b) is much simpler, so we have decided to consider it separately.

\medskip

\noindent
{\it Case (a): $n=4$.\ \ }  
In this case $\PP$ is a tetrahedron with triangular faces 
$T_1, T_2, T_3$, and $T_4$. Without loss of generality we may assume that  
$\BB_2^3$ is tangential to the faces $T_1, T_2, T_3$, and $T_4$. Then the 
well-known inequality between the harmonic and arithmetic means yields that
$$
  1=\sum_{i=1}^{4}\frac{\frac{1}{3}{\rm area}(T_i)}{{\rm vol}(\PP)}\ge
  \sum_{i=1}^{4}\frac{1}{|p_i|+1}\ge\frac{4^2}{\sum_{i=1}^{4}(|p_i|+1)}.
$$
This implies in a straightforward way that
$$
  \sum_{i=1}^{4}|p_i|\ge 12>6\sqrt{3},
$$
finishing the proof of this case.

\medskip

For the next two cases we will need the following notation. Fix $i\leq n$. 
Let $C_i$ denote the (closed) spherical cap of $\S^2$ with spherical radius 
$R_i$ which is the union of points $x\in \S^2$ such that the open line segment 
connecting $x$ and $p_i$ is disjoint from $\BB_2^3$. In other words, $C_i$ is 
the spherical cap with the center $p_i/|p_i|$ and the spherical radius $R_i$, 
satisfying $|p_i|=\frac{1}{\cos R_i}$. By $b_i$ we denote the spherical 
area of $C_i$. Then  $b_i = 2\pi (1-\cos R_i)$. 

\medskip

\noindent
{\it Case (b): $n\geq 8$.\ \ } 
Since $\PP$  contains $\BB_2^3$, we have 
$$
   \S^2 \subset \bigcup _{i=1}^{n} C_i.
$$  
Comparing the areas, we observe 
$$
  4 \pi \leq \sum _{i=1}^{n} b_i = \sum _{i=1}^{n}  2\pi \left(1-\cos R_i \r), 
$$
which implies 
$$
  \sum _{i=1}^{n} \cos R_i \leq n-2. 
$$
Applying again the inequality between the harmonic and arithmetic means, we 
obtain
$$
   \sum _{i=1}^{n}  |p_i|= \sum _{i=1}^{n}  \frac{1}{\cos R_i} \geq 
   \frac{n^2}{\sum _{i=1}^{n}  \cos R_i} \geq \frac{n^2}{n-2} \geq 
    \frac{64}{6} > 6\sqrt{3}. 
$$

\medskip

\noindent
{\it Case (c): $5\leq n\leq 7$.\ \ } 
Let $\PP^{\circ}$ denote the polar of $\PP$. Given $i\leq n$,  
let $F_i$ denote the central projection of the face of $\PP^{\circ}$ that 
corresponds to the vertex $p_i$ of $\PP$ from the center $0$ onto the boundary of $\BB_2^3$,  
i.e. onto the unit sphere $\S^2$ centered at $0$. 
Obviously, $F_i$ is a spherically convex polygon of $\S^2$ and $F_i \subset C_i$. 
Let $n_i$ denote the number of sides of $F_i$ and let $a_i$ stand for the 
spherical area of $F_i$. Note that the area of the sphere is equal to the sum of 
areas of $F_i$'s, that is $\sum_{i=1}^{n} a_i =4\pi$. 
As $10<6\sqrt{3}=10.3923...<11$, therefore without loss of generality we may assume that 
there is no $i$ for which $|p_i|=\frac{1}{\cos R_i}\ge 11-3=8$, in other words we assume 
that $0<R_i<\arccos\frac{1}{8}=1.4454...<\frac{\pi}{2}$ 
for all $i\le n$. 
Note that this immediately implies that 
$0<a_i<b_i = 2\pi (1-\cos R_i) < \frac{7\pi}{4}< 5.5$ 
for all $1\le i\le n$.

It is well-known that if $C\subset \S^2$ is a (closed) spherical cap of 
radius less than $\frac{\pi}{2}$, then the spherical area of a spherically convex 
polygon with 
at most $s\ge 3$ sides lying in $C$ is maximal for the regular spherically convex 
polygon with $s$ sides inscribed in $C$. (This can be easily obtained with the help 
of the Lexell-circle (see \cite{F}).)  
It is also well-known that if $F_i^*$ denotes a regular spherically convex polygon with 
$n_i$ sides and of spherical area $a_i$, and  if $R_i^*$ denotes the circumradius of 
$F_i^*$, then 
$\frac{1}{\cos R_i^*}=\tan \frac{\pi}{n_i}\tan\big(\frac{a_i+(n_i-2)\pi}{2n_i}\big)$. 
Thus, for every $i\leq n$ we have 
$$
   |p_i|=\frac{1}{\cos R_i}\ge \tan \frac{\pi}{n_i} 
   \tan\left(\frac{a_i+(n_i-2)\pi}{2n_i}\right) .
$$
Here $3\le n_i\le n-1\le 6$ and $0<a_i<\frac{7\pi}{4}$
for all $1\le i\le n$. 
 
Now, it is natural to consider the function 
$f(x,y)=\tan \frac{\pi}{y}\tan\big(\frac{x+(y-2)\pi}{2y}\big)$  
defined on $\{(x,y)\ |\ 0< x < 2\pi ,\  3\le y\}$. 
As in 2-dimensional case we are going to use the Jensen inequality.  
But, unfortunately, it turns out that $f$ is convex only on a proper 
subset of its domain, see Lemma~\ref{func}. 
Without loss of generality we may assume that $m$ is chosen such that 
$0<a_i<0.4$ for all $i\le m$ and $0.4\le a_i <5.5$ for all $m+1\le i\le n$. 
Since $\sum _{i=1}^n a_i = 4\pi$, one has $m<n-1$. 
By Lemma~\ref{func} (iv) and by the Jensen inequality, we obtain 
$$
   \sum_{i=1}^{n} |p_i|\ge \sum_{i=1}^{m} |p_i|+\sum_{i=m+1}^{n} f(a_i,\  n_i)
$$
$$
  \ge m+(n-m)\ f\left(\frac{1}{n-m}
   \ \sum_{i=m+1}^n a_i,\  \frac{1}{n-m}\ \sum_{i=m+1}^n n_i \right) 
$$
(here by $\sum_i^0$ we mean $0$). Since $\sum_{i=1}^{n} a_i =4\pi $, we have 
$\sum_{i=m+1}^n a_i>4\pi -0.4m$. By Euler's theorem on the edge graph of 
$\PP^{\circ}$ we also have that $\sum_{i=1}^{n} n_i\le 6n-12$ and therefore 
$\sum_{i=m+1}^{n} n_i\le (6n-12)-3m$. 
Thus, applying Lemma~\ref{func} (i) and (ii), we observe
$$
  \sum_{i=1}^{n} |p_i| \ge m+(n-m)f\left(\frac{4\pi-0.4m}{n-m},\  
  \frac{(6n-12)-3m}{n-m}\right) =: g(m,n). 
$$

First we show that $g(m,n)\ge 6\sqrt{3}=10.3923...$ for every 
$(m, n)$ with $6 \leq n \leq 7$ and $0\leq m <n-1$.

\smallskip

\noindent{\it Subcase $n=7$:} 
$$g(0,7)=10.9168..., \quad g(1,7)=10.8422..., \quad g(2,7)=10.8426...,$$
$$g(3,7)=11.0201..., \quad g(4,7)=11.7828..., \quad g(5,7)=18.3370... .$$

\smallskip

\noindent
{\it Subcase $n=6$:} 
$$g(0,6)=6\sqrt{3}=10.3923...,  \quad g(1,6)=10.4034..., \quad g(2,6)=10.6206..., $$
$$g(3,6)=11.5561..., \quad g(4,6)=21.2948... .$$

\smallskip

\noindent
{\it Subcase $n=5$:}
First note that 
$$6\sqrt{3}< g(1,5)=10.6302...< g(2,5)=11.8680... < g(3,5)=28.1356... .$$ 
Unfortunately, $g(0, 5) < 6\sqrt{3}$, so we treat the case $n=5$ slightly 
differently (in fact the proof is easier than the proof of the case 
$6\leq n\leq 7$, since we will use convexity of a function of one variable). 

 In this case $\PP$ has only 5 vertices, so it is either a double tetrahedron 
or a cone over a quadrilateral. As the later one can be thought of as a limiting 
case of double tetrahedra, we can assume that the edge graph of $\PP$ has two 
vertices, say $p_1$ and $p_2$, of degree three and three vertices, say $p_3, p_4$, 
and $p_5$, of degree four. Thus $n_1=n_2=3$ and $n_3=n_4=n_5=4$. Therefore  
$$
  \sum_{i=1}^{5}|p_i|\ge \sum_{i=1}^{5} f(a_i,\  n_i) = 
  \sum_{i=1}^{2} f(a_i,\  3) + \sum_{i=3}^{5} f(a_i,\  4) .  
$$
By Lemma~\ref{func} (iii) and by the Jensen inequality, we get  
$$ 
 \sum_{i=1}^{5} |p_i| \ge 2\ f\left( \frac{a_1+a_2}{2},\ 3 \right) + 3 \ 
  f\left( \frac{a_3+a_4+a_5}{3},\  4\right)  
$$ 
$$ 
  = 2\ f(a,\ 3)+3 \ f\left( \frac{4\pi-2a}{3},\ 4\right) 
$$
$$ 
  = 2\sqrt{3}\ \tan\left(
  \frac{a+\pi}{6}\right) + 3 \ \tan\left( \frac{5\pi-a}{12}\right) =: h(a),
$$
where $0 \leq a =\frac{a_1+a_2}{2}< 5.5$. 
Finally, it is easy to show that the minimum value of $h(a)$ over the closed interval 
$0\le a\le 5.5$ is (equal to $10.5618...$ and therefore is) strictly larger than 
$6\sqrt{3}=10.3923...$, 
completing the proof of the first part of the theorem. 

\smallskip

$(ii)$ First, observe that $(i)$, John's Theorem, and Lemma~\ref{distest} imply that
$$
  4=\frac{4\sqrt{2}}{\sqrt{2}}\le\frac{{\rm vein}(\BB_2^2)}{d_{\KK} }\le{\rm vein}(\KK)
$$
   and  
$$ 
  6=\frac{6\sqrt{3}}{\sqrt{3}}\le\frac{{\rm vein}(\BB_2^3)}{d_{\LL} }
  \le{\rm vein}(\LL) \le{d_{\LL} }\cdot{\rm vein}(\BB_2^3)\le 18.
$$
Second, Corollary~\ref{ill} shows that indeed ${\rm vein}(\KK)\le 6$,  
finishing the proof. 
\kkk

\medskip

\noindent 
{\bf Remark. } 
Note that the proof of Case (a) works in higher dimensions as well. 
Namely, if $\PP$ is a simplex containing the Euclidean ball $\BB_2^d$ 
and the $p_i$'s denote the vertices of $\PP$, then 
$$
   \sum _{i=1}^{n+1} |p_i| \geq d (d+1)
$$ 
 with  equality only for regular simplices circumscribed $\BB_2^d$. 

\medskip

\section{The vertex index in the high dimensional case}

In this section we deal with the high dimensional case. First, we compute precisely 
${\rm vein}(\BB _1^d)$. Then we provide a lower and an upper estimates in the general 
case.

In fact, the estimate for  ${\rm vein}(\BB _1^d)$ follows now from the more general 
fact, namely ${\rm vein}(\KK) \geq 2d$ for every $0$-symmetric $\KK$ in $\R^d$, 
proved in \cite{GL}. However the proof of this fact is very non-trivial and quite long, 
so we have decided to present a simple direct proof for the case 
$\KK = \BB _1^d$.  

\begin{prop} \label{octah}
For every $d\ge 2$ one has
$$
   {\rm vein}(\BB_1^d) = 2d. 
$$
\end{prop}

\medskip

\proof
The estimate ${\rm vein}(\BB_1^d) \leq 2d$ is trivial.

Now assume that $\{p_i\} _{i=1}^N$ be such that $p_i = \{p_{i j}\}_{j=1}^d \in \R^d$ for 
every $i\leq N$ and $\BB_1^d \subset \conv  \{p_i\} _{i=1}^N$. Then for every $k\leq d$ 
we have that $e_k$ and $-e_k$ are convex combinations of $p_i$'s, that is, there are 
$\{\alpha _{k i}\} _{i=1}^N$ and $\{\beta _{k i}\} _{i=1}^N$ such that $\alpha _{k i}\geq 0$, 
$\beta _{k i} \geq 0$, $i\leq N$, and 
$$
   \sum _{i=1}^N \alpha _{k i} = \sum _{i=1}^N \beta _{k i} = 1, \quad e_k = 
   \sum _{i=1}^N \alpha _{k i} p _i, \quad -e_k = \sum _{i=1}^N \beta _{k i} p_i. 
$$
It implies for every $k$ 
$$
    1 = \sum _{i=1}^N \alpha _{k i} p _{i k} \leq \max _{i\leq N} p_{i k}
$$
and 
$$
   - 1 = \sum _{i=1}^N \beta _{k i} p _{i k} \geq \min _{i\leq N} p_{i k}.
$$
Therefore
$$
    \sum _{i=1}^N \| p _{i } \| _1 = \sum _{i=1}^N \sum _{k=1}^d | p_{i k} | \geq 
    \sum _{k=1}^d \left( \max _{i\leq N} p_{i k} - \min _{i\leq N} p_{i k} \right) 
    \geq 2 d, 
$$
proving the lower estimate ${\rm vein}(\BB_1^d) \geq 2d$.
\kkk

\subsection{A lower bound}

In this section, we provide a lower estimate for ${\rm vein}(\KK)$ in terms 
of {\it outer volume ratio} of $\KK$. As the example of 
the Euclidean ball shows,  
our estimate can be asymptotically sharp. 

\begin{theorem} \label{ball}
There is an absolute constant $c>0$ such that for every $d\ge 2$ 
and every $0$-symmetric convex body $\KK$ in $\R^d$ one has
$$
  \frac{  d^{3/2}}{\sqrt{2 \pi e}\ \mbox{\rm ovr} (\KK) } \leq  
  \frac{ \ d}{ \left( \vol (\BB_2^d)\r)^{1/d} \ \mbox{\rm ovr} (\KK)} 
  \leq {\rm vein}(\KK) .
$$
\end{theorem}

\proof
Recall that ${\rm vein}(\KK)$ is an affine invariant, i.e. 
${\rm vein}(\KK) = {\rm vein}(T \KK)$ for every invertible linear 
operator $T : \R^d \to \R^d$. Thus, without lost of generality 
we can assume that $\BB_2^d$ is the ellipsoid of minimal volume
for $\KK$. In particular, $\KK \subset \BB_2^d$, so 
$|\cdot |\leq \no _{\KK}$. 

Let $\{p_i\}_1^N\in \R^d$ be such that 
$\KK \subset \conv \{p_i\}_1^N$. Clearly $N\geq d+1$. Denote 
$$
  \LL := \absconv \{ p_i \}_1^N. 
$$
Then 
$$
 \LL^{\circ} = \left\{ x \, \, \mid \, \, |\la x, p_i\ra | \leq 1 \, \, 
  \mbox{ for every } \, \, i \leq N \r\}. 
$$
By Theorem~2 of \cite{BP}, we observe 
$$
  \mbox{\rm vol} \left( \LL^{\circ} \r) \geq 
  \left( \frac{d}{\sum _1^N |p_i|} \r)^d . 
$$ 
Since, by Santal\'o inequality $\mbox{\rm vol} \left( \LL \r) 
\mbox{\rm vol} \left( \LL^{\circ} \r) \leq \left(\mbox{\rm vol} 
\left( \BB_2^d  \r) \r)^2$ and since $\KK \subset \LL$, we obtain 
$$
   \mbox{\rm vol} \left( \KK \r) \leq \mbox{\rm vol} \left( \LL \r) 
  \leq   \frac{ \left(\mbox{\rm vol} \left( \BB_2^d  \r) \r)^2}{ \mbox{\rm vol}
  \left( \LL^{\circ} \r)}  \leq \left(\mbox{\rm vol} \left( \BB_2^d  
  \r) \r)^2 \left( \frac{1}{d}  \ \sum _1^N |p_i|  \r) ^d . 
$$ 
Finally, since $\BB_2^d$ is the minimal volume ellipsoid for $\KK$ and 
since $\no_ {\KK}\geq |\cdot|$, we have 
$$
   \frac{1}{\mbox{\rm ovr} (\KK)} = \left( \frac{ \mbox{\rm vol} \left(  \KK \r) }{ 
   \mbox{\rm vol} \left( \BB_2^d \r)} \r)^{1/d} \leq \left( \mbox{\rm vol} 
   \left(  \BB_2^d \r) \r)^{1/d} \ \frac{1}{d}  \ \sum _1^N \|p_i\|_ {\KK} ,  
$$  
which implies the desired result.
\kkk

\medskip

We have the following immediate corollary of Theorem~\ref{ball}.

\begin{sled} \label{ballest} 
For every $d\ge 2$ one has
$$
  \frac{  d^{3/2}}{\sqrt{2 \pi e}}  \leq {\rm vein}(\BB_2^d) \leq 2 \ d^{3/2}, 
 \quad \quad \quad 
  \frac{  d^{3/2}}{\pi e} \leq {\rm vein}(\BB_{\infty}^d) \leq C\ d^{3/2},  
$$
where $C=2$ if $d=2^m$ for some integer $m$ and $C=2/(\sqrt{2}-1)$ in general. 
\end{sled}

\proof
The lower estimates here follow from Theorem~\ref{ball} and computation of volumes.
Indeed, as we noticed above, 
$$
  \vol (\BB_2^d) \leq \left(\frac{2 \pi e}{d}\r)^{d/2}
$$
and, therefore, 
$$
  \mbox{\rm ovr} (\BB _{\infty}^d) \leq \left(\frac{\vol \left( \sqrt{d} \ 
  \BB_2^d\r)}{\vol \left( \BB_{\infty}^d\r)} \r) ^{1/d} \leq \frac{\sqrt{2 \pi e}}{2}. 
$$ 

The upper estimates 
follow from Proposition~\ref{octah} and Lemma~\ref{distest}, 
since $d (\BB_2^d, \BB_1^d) = \sqrt{d}$ 
and, by Theorem~\ref{cubeoc}, $d (\BB_{\infty}^d, \BB_1^d)\leq (C/2) \sqrt{d}$. 
\kkk

\medskip

\subsection{An upper bound}

Let $u$, $v \in \R^d$. 
As usual $Id  :  \R^d \to \R^d$ denotes the identity operator and $u\otimes v$ denotes the 
operator from $\R^d$ to $\R^d$, defined by $(u\otimes v) (x) = \la u, x \ra v$ for every 
$x\in \R^d$. In \cite{R1, R2}, M.~Rudelson proved the following theorem (see Corollary~4.3 
of \cite{R1} and Theorem~1.1 with Remark~4.1 of \cite{R2}).

\begin{theorem}\label{rudel}
For every $0$-symmetric convex body $\KK$ in $\R^d$ and every 
$\eps \in (0, 1]$ there exists a $0$-symmetric convex body 
$\LL$ in $\R^d$ such that $d(\KK, \LL)\leq 1+\eps$ and $\BB_2^d$ is 
the minimal volume ellipsoid containing $\LL$, and
$$
   Id = \sum_{i=1}^{M} c_i u_i\otimes u_i, 
$$
where $c_1$, ..., $c_M$ are positive numbers, $u_1$, ..., $u_M$ are contact points 
of $\LL$ and $\BB_2^d$ (that is $\|u_i\| _{\LL}=|u_i|=1$), and 
$$
   M\leq C \ \eps ^{-2}\ d  \ \ln(2d),  
$$
with an absolute constant $C$.
\end{theorem}

\noindent
{\bf Remark.} It is a standard observation (cf. \cite{Ball}, \cite{T}) 
that under the conditions of Theorem~\ref{rudel} one has
$$
   \PP \subset \LL \subset \BB_2^d \subset \sqrt{d} \ \LL, 
$$
for $\PP = \absconv \{u_i\}_{i\leq M}  $. Indeed, $\PP \subset \LL$ by the 
convexity and the symmetry of $\LL$, and for every $x\in \R^d$ we have 
$$
  x = Id \ x =  \sum_{i=1}^{M} c_i \la u_i, x\ra  u_i, 
$$
so 
$$
  |x|^2 = \la x, x\ra  =  \sum_{i=1}^{M} c_i \la u_i, x\ra ^2 \leq 
     \max _{i\leq M }\la u_i, x\ra ^2  \sum_{i=1}^{M} c_i 
   = \|x \| ^2_{\PP ^{\circ}}  \  \sum_{i=1}^{M} c_i .
$$
Since 
$$
  d = \mbox{trace}\  Id = \mbox{trace}   \sum_{i=1}^{M} c_i u_i\otimes u_i = 
   \sum_{i=1}^{M} c_i \la u_i, u_i \ra =  \sum_{i=1}^{M} c_i ,  
$$
we obtain $|x|\leq \sqrt{d} \  \|x \| _{\PP ^{\circ}}$, which means 
$\PP ^{\circ}\sqrt{d} \subset B_2^d$. 
By duality we have $B_2^d \subset \sqrt{d} \ \PP$.
Therefore, $d(\KK, \PP) \leq d(\KK, \LL)\ d(\LL, \PP)\leq (1+\eps) \sqrt{d}$, 
and, hence, we have the following immediate consequence of Theorem~\ref{rudel}.

\begin{sled} \label{rud}
For every $0$-symmetric convex body $\KK$ in $\R^d$ and every $\eps \in (0, 1]$ there 
exists a $0$-symmetric convex polytope $\PP$ in $\R^d$ with $M$ vertices such that 
$d(\KK, \PP)\leq (1+\eps)\sqrt{d}$ and 
$$
   M\leq C \ \eps ^{-2}\ d  \ \ln(2d),  
$$
where $C$ is an absolute constant.
\end{sled}

This Corollary implies the general upper estimate for ${\rm vein}(\KK)$.

\begin{theorem} \label{ru}
For every centrally symmetric convex body $\KK$ in $\R^d$ one has 
$$
   {\rm vein} (\KK) \leq C\ d^{3/2}  \ \ln(2d),  
$$
where $C$ is an absolute constant.
\end{theorem}

\proof Let $\PP$ be a polytope given by Corollary~\ref{rud} applied 
to $\KK$ with $\eps =1$. Then  $d(\KK, \PP)\leq 2 \sqrt{d}$. Clearly, 
${\rm vein}(\PP) \leq M$ (just take the $p_i$'s in the definition of ${\rm vein}(\cdot)$ to 
be vertices of $\PP$). Thus, by Lemma~\ref{distest}  we obtain 
${\rm vein}(\KK) \leq 2  M \sqrt{d}$, which completes the proof. 
\kkk

\vspace{1cm}

\medskip

\noindent
K\'aroly Bezdek,
Department of Mathematics and Statistics,
2500 University drive N.W.,
University of Calgary, AB, Canada, T2N 1N4.
\newline
{\sf e-mail: bezdek@math.ucalgary.ca}
 
\smallskip

\noindent
A.E. Litvak,
Department  of Mathematical and Statistical Sciences,
University of Alberta, Edmonton, AB, Canada, T6G 2G1.
\newline
{\sf e-mail: alexandr@math.ualberta.ca}


\begin{thebibliography}{GGM}


\bibitem[A]{A} E.~Asplund, 
{\em  Comparision between plane symmetric convex bodies and 
parallelograms}, Math. Scand. 8 (1960), 171--180. 
%
%
\bibitem[Ba]{Ball} K.~Ball, 
{\em  Flavors of geometry} in 
{\em An elementary introduction to modern convex geometry}, 
Levy, Silvio (ed.), Cambridge: 
Cambridge University Press. Math. Sci. Res. Inst. Publ. 31, 1--58 (1997). 
%
%
\bibitem[BaP]{BP}  K.~Ball, A.~Pajor, 
{\em Convex bodies with few faces}, 
Proc. Am. Math. Soc. 110 (1990),  225--231. 
%
%
\bibitem[Be1]{Be92} K.~Bezdek, 
{\em Research problem 46}, Period. Math. Hungar.  24 (1992), 119--121.
%
%
\bibitem[Be2]{Be93} K.~Bezdek, {\em Hadwiger-Levi's covering problem revisited}, 
New Trends in Discrete and Computational Geometry (ed.: J. Pach), Springer-Verlag, 1993.
%
%
\bibitem[Be3]{Be05} K.~Bezdek, {\em The illumination conjecture and its extensions}, 
Period. Math. Hungar. 53 (2006), 59--69.  
%
%
\bibitem[BeBK]{Be06} K.~Bezdek, K.~B\"or\"oczky, Gy.~Kiss, 
{\em On the successive illumination parameters of convex bodies}, 
Period. Math. Hungar. 53 (2006), 71--82.
%
%
\bibitem[F]{F} L.~Fejes T\'oth, {\em Lagerungen in der Ebene, auf der Kugel und im Raum},
Springer Verlag, Berlin-G\"ottingen-Heidelberg, 1953.
%
%
\bibitem[GL]{GL}
 E.~D.~Gluskin, A.~E.~Litvak, 
{\em Asymmetry of convex polytopes and 
vertex index of symmetric convex bodies}, preprint. 
%
%
\bibitem[GMP]{GMP}
 Y.~Gordon, M.~Meyer, A.~Pajor, 
{\em Ratios of volumes and factorization through $\ell _{\infty}$}, 
Ill. J. Math. 40 (1996),  91--107. 
%
%
\bibitem[GKM]{GKM}
V.~I.~Gurari, M.~I.~Kadec, V.~I.~Macaev,  
{\em Distances between finite-dimensional analogs of the $L_{p}$-spaces}, 
(Russian),  Mat. Sb. (N.S.),  70 (112),  1966, 481--489. 
%
%
\bibitem[J]{J} F.~John, 
{\em Extremum problems with inequalities as subsidiary conditions},
Studies and Essays Presented to R. Courant on his 60th Birthday,
January 8, 1948, 187--204. Interscience Publishers, Inc., New York,
N. Y., 1948.
%
%
\bibitem[MS]{MS} H.~Martini, V.~Soltan, {\em Combinatorial problems on the 
illumination of convex bodies}, Aequationes Math. 57 (1999), 121--152.
%
%
\bibitem[Pi]{Pisier} G.~Pisier,
{\em The Volume of Convex Bodies and Banach Space Geometry}, 
Cambridge University Press 1989.
%
%
\bibitem[R1]{R1}
M.~Rudelson, 
{\em Random vectors in the isotropic position},  
J. Funct. Anal. 164 (1999),  60--72. 
%
%
\bibitem[R2]{R2} 
M.~Rudelson,  
{\em Contact points of convex bodies},  
Israel J. Math. 101 (1997), 93--124. 
%
%
\bibitem[Sw]{S}K.~J.~Swanepoel, 
{\em Quantitative illumination of convex bodies and vertex degrees of 
geometric Steiner minimal trees}, Mathematika 52 (2005), 47--52. 
%
%
\bibitem[T]{T} N.~Tomczak-Jaegermann,
{\em Banach-Mazur distances and finite-dimensional operator ideals},
Pitman Monographs and Surveys in Pure and Applied Mathematics, 38.
Longman Scientific \& Technical, Harlow; copublished in the
United States with John Wiley \& Sons, Inc., New York, 1989.
%
%
\end{thebibliography}
\end{document}